\documentclass[sts]{imsart}

\RequirePackage[OT1]{fontenc}
\usepackage{amsthm,amsmath,amssymb,natbib,xspace,graphicx}
\RequirePackage[colorlinks,citecolor=blue,urlcolor=blue]{hyperref}

\startlocaldefs
\numberwithin{equation}{section}
\theoremstyle{plain}

\DeclareSymbolFont{GreekLetters}{OML}{cmr}{m}{it} %Provide missing letters
\DeclareSymbolFont{UpSfGreekLetters}{U}{cmss}{m}{n} %Provide missing letters
\DeclareMathSymbol{\varrho}{\mathalpha}{GreekLetters}{"25}
\DeclareMathSymbol{\UpSfLambda}{\mathalpha}{UpSfGreekLetters}{"03}
\DeclareMathSymbol{\UpSfSigma}{\mathalpha}{UpSfGreekLetters}{"06}
\providecommand{\mathbold}{\boldsymbol}
\newcommand{\bvec}[1]{\mathbold{#1}}

\newlength{\overwdth}

\def\abs#1{\ensuremath{\left \lvert #1 \right \rvert}}

\newcommand{\norm}[2][{}]{\ensuremath{\left \lVert #2 \right \rVert}_{#1}}

\newcommand{\IIDsim}{\overset{\textup{IID}}{\sim}}

\DeclareMathOperator{\Prob}{\mathbb{P}}
\DeclareMathOperator{\Ex}{\mathbb{E}}

\DeclareMathOperator{\diag}{diag}

\DeclareMathOperator*{\argmin}{argmin}

\DeclareMathOperator{\Order}{{\mathcal O}}

\newcommand{\vone}{\bvec{1}}

\newcommand{\va}{\bvec{a}}
\newcommand{\vA}{\bvec{A}}
\newcommand{\vb}{\bvec{b}}

\newcommand{\vc}{\bvec{c}}

\newcommand{\vf}{\bvec{f}}

\newcommand{\vk}{\bvec{k}}
\newcommand{\vK}{\bvec{K}}

\newcommand{\vell}{\bvec{\ell}}

\newcommand{\vv}{\bvec{v}}
\newcommand{\vV}{\bvec{V}}

\newcommand{\vx}{\bvec{x}}

\newcommand{\vtheta}{\bvec{\theta}}

\newcommand{\tf}{\widetilde{f}}

\newcommand{\tvf}{\tilde{\vf}}

\newcommand{\naturals}{\mathbb{N}}

\newcommand{\cg}{\mathcal{G}}
\newcommand{\ch}{\mathcal{H}}

\newcommand{\cn}{\mathcal{N}}
\newcommand{\cp}{\mathcal{P}}

\newcommand{\calH}{\mathcal{H}}
\newcommand{\vC}{\boldsymbol{C}}
\newcommand{\calGP}{\cg\!\cp}
\DeclareMathOperator{\Var}{\mathbb{V}}
\defcitealias{BriEtal18a}{PI}% Define a citation alias
\newcommand{\BOGOS}{\citetalias{BriEtal18a}}% Use shorthand
\defcitealias{RatHic19a}{JH}% Define a citation alias
\newcommand{\JH}{\citetalias{RatHic19a}}% Use shorthand
\newcommand{\vLambda}{\boldsymbol{\Lambda}}
\endlocaldefs

\begin{document}

\begin{frontmatter}
\title{Comment on ``Probabilistic Integration: A Role in Statistical Computation?''
\thanksref{T1}}
\runtitle{Comment on ``Probabilistic Integration \ldots''}
\thankstext{T1}{This work is supported in part by NSF-DMS-1522687.}

\begin{aug}
\author{\fnms{Fred J.} \snm{Hickernell}\ead[label=e1]{hickernell@iit.edu} \ead[label=u1,url]{iit.edu/$\sim$hickernell}}
\and
\author{\fnms{R.} \snm{Jagadeeswaran}\ead[label=e2]{jrathin1@hawk.iit.edu}}

\runauthor{Fred J. Hickernell and R. Jagadeeswaran}

\affiliation{Illinois Institute of Technology}

\address{RE 208, 10 W.\ 32$^{\text{nd}}$ St., Chicago, IL 60616
\printead{e1,e2,u1}.}

\end{aug}

\begin{abstract} 
Probabilistic integration provides a criterion for stopping a simulation when a specified error tolerance is satisfied with high confidence.  We comment on some of the modeling assumptions and implementation issues involved in designing an automatic Bayesian cubature.
\end{abstract}


\end{frontmatter}

\section{When to Stop?}

In highlighting the possibilities of probabilistic integration, the authors of \cite{BriEtal18a}, henceforth abbreviated as \BOGOS{}, have suggested a useful stopping criterion for cubature.  Numerical analysis provides an upper bound on the cubature error expressed as a product of the roughness of the integrand and the quality of our sampling scheme.  For example, \BOGOS{} (5) quotes the error bound
\begin{equation} \label{HJErrBd}
    \bigl \lvert \hat{\Pi}[f] - \Pi[f] \rvert \le \lVert f \rVert_{\calH} \lVert \mu(\hat{\pi}) - \mu(\pi) \rVert_{\calH},
\end{equation}
where 
\begin{itemize}
    \item the integrand, $f$, lies in a Hilbert space, $\ch$,
    \item $\Pi[f]$ denotes the desired integral of $f$ defined in terms of the probability measure $\pi$, and
    \item $\hat{\Pi}[f]$ denotes a cubature defined in terms of the discrete measure $\hat{\pi}$.
\end{itemize}
The \emph{discrepancy} between $\pi$ and $\hat{\pi}$ is defined as $\lVert \mu(\hat{\pi}) - \mu(\pi) \rVert_{\calH}$. As the sample size, $n$, increases, a well chosen sequence of discrete measures causes the discrepancy to tend to zero.  

But, even if $\lVert \mu(\hat{\pi}) - \mu(\pi) \rVert_{\calH}$ can be computed efficiently, one typically does not have a good estimate or bound on $\lVert f \rVert_{\calH}$.  Therefore, it is impractical to use  \eqref{HJErrBd} to determine an $n$ satisfying the error criterion
\begin{equation} \label{HJErrCrit}
    \bigl \lvert \hat{\Pi}[f] - \Pi[f] \bigr \rvert \le \varepsilon,
\end{equation}
where $\varepsilon$ is the user-specified absolute error tolerance.  

We believe that the practitioner would like an automatic cubature, i.e., an algorithm with a stopping criterion that guarantees \eqref{HJErrCrit} (with high probability).  Probabilistic integration, and in particular Bayesian cubature, as espoused in \BOGOS{}, fulfills that wish.  

Bayesian cubature, as explained in \BOGOS{}, assumes that the integrand, $f$, may be modeled by a Gaussian stochastic process, $g \sim \calGP(0,c)$, conditioned on $g$ having the same values as $f$ at the cubature nodes or states, $\{\vx_i\}_{i=1}^n$.  Thus, $\hat{\Pi}[g] = \hat{\Pi}[f]$.  Furthermore, Bayesian cubature is designed to satisfy $\hat{\Pi}[g] = \Ex_n[\Pi[g]]$.  Here, $c$ is the covariance function (or kernel) for $g$.  The definition of $g$ allows us to construct credible intervals for the cubature error via Proposition 1 in  \BOGOS{}, namely,
\begin{gather}
\label{credInt}
    \Prob\Bigl[\bigl \lvert \hat{\Pi}[f] - \Pi[g] \bigr \rvert \le 2.58 \sqrt{\Var_n[\Pi[g]]}  \Bigr] = 99\%, \\
    \label{Vnform}
    \Var_n[\Pi[g]] = \Pi\Pi[c(\cdot,\cdot)] - \Pi[\vc(\cdot,X)] \vC^{-1} \Pi[\vc(X,\cdot)].
\end{gather}
If the observed integrand, $f$, lies in the $99\%$ middle of the sample space for $g$, and not in the $1\%$ extreme, then increasing $n$ until $2.58 \sqrt{\Var_n[\Pi[g]]}$ is no greater than $\varepsilon$ ensures that  \eqref{HJErrCrit} holds with $99\%$ probability.

There are some practical obstacles to implementing this elegant recipe. 

\begin{itemize} 

\item How does one choose the covariance function $c$?  While one may always choose the sample space large enough to include $f$, our use of the credible interval as a stopping criterion assumes that $f$ is not in the tails of the distribution $\calGP(0,c)$.  We discuss this question in the next section.  

\item The computational cost of computing $\Var_n[\Pi[g]]$ involves matrix inversion, which requires $\mathcal{O}(n^3)$ operations in general.  This typically takes much more time than the $\mathcal{O}(n)$ operations required to compute the cubature, $\hat{\Pi}[f]$, unless obtaining an integrand value is quite time-consuming. We discuss how to circumvent this problem by matching covariance functions and cubature nodes in Section \ref{sec:Match}.

\end{itemize}
There are commonalities and differences in the deterministic and Bayesian approaches to numerical integration.  We discuss some of these in Section \ref{sec:ProbDet}.

%%%%%%%%%%%%%%%%%%%%%%%%%%%%%%%%%%%
\section{Which Gaussian Process?} \label{sec:WhichGauss}
As mentioned above, using a credible interval as a stopping criterion requires a careful choice of the covariance function, $c$.  The width of the credible interval in \eqref{credInt} depends on $\Var_n[\Pi[g]]$ given by \eqref{Vnform}.  At first glance, nothing in \eqref{Vnform} depends on the integrand data, $\vf = \bigl( f(\vx_i) \bigr)_{i=1}^n$, although our intuition tells us that it should.  The credible interval for the integral of $47f$ should be $47$ times as wide as the credible interval for the integral of $f$.  

When constructing the confidence interval for the mean of a scalar random variable, $Y$, from independent and identically distributed (IID) data, one must estimate the variance of $Y$ by the sample variance. Analogously, when constructing the credible interval in \eqref{credInt} for the integral (mean) of a function, one must estimate the vertical scale factor inherent in the covariance function $c$.

We have recently explored Bayesian cubature as the basis for automatically selecting $n$ to satisfy the error criterion \eqref{HJErrCrit}  in \cite{RatHic19a}, henceforth abbreviated as \JH{}.  As in Proposition 2 of \BOGOS{},  \JH{} chooses the covariance function to take the form $c(\vx, \vx') =  \lambda c_0(\vx, \vx'; \vtheta)$, where $\lambda$ is the vertical scale factor, and the parameter $\vtheta$ determines the smoothness and other properties of the covariance function.  An example of $c_0$ is  the following (\JH{} (36)):
\begin{multline}
\label{the_kernel_eqn_bernoulli}
c_0(\vx, \vx';\vtheta) =
\prod_{l=1}^d \biggl[
1 - (-1)^{r} \gamma B_{2r}( |{x_l-x'_l}| ) \biggr], \\  
\forall \vx,\vx' \in [0,1]^d, \  \vtheta = (r,\gamma), \ r \in \naturals, \ \gamma > 0,
\end{multline}
where $B_{2r}$ is the Bernoulli polynomial of degree $2r$.  The smoothness of the covariance function increases with $r$.  Covariance functions of this form appear in  \cite{Hic96a,DicEtal14a}. Bernoulli polynomials are described in Chapter 24 of \cite{OlvEtal10a}.

To increase the possibility that our integrand $f$ lies in the middle of the sample space, we also allow the Gaussian process $g$ to have an arbitrary mean, $m$, so $g \sim \calGP(m, \lambda c_0)$.  One may imagine the situation where $f$ represents an option payoff.  Then, $f$ is non-negative and its mean is non-negative.  Assuming an improper prior on $(m, \lambda)$, the posterior marginal for $\Pi[g]$ is  a Student-t distribution with $n-1$ degrees of freedom and with mean and variance both depending on the integrand data, $\vf$ (\JH{} (15--16)):
\begin{subequations} \label{JHFullBayes}
\begin{align}
\label{eqn:BC}
    \hat{\Pi}[f] & =  \Ex_n[\Pi[g]] \\
    \nonumber
    & =
\left(
\frac{  (1 - \vone^T  \vC^{-1}_0\Pi[\vc_0(X,\cdot) ] \vone^T }{ \vone^T  \vC^{-1}_0 \vone}   +  \Pi[\vc_0(\cdot, X)]
\right)  \vC^{-1}_0 \vf, \\
\label{eqn:errMLE}
\Var_n[\Pi[g]] & = \frac{1}{n-1}
 \vf^T \left(  \vC^{-1}_0 - 
\frac{  \vC^{-1}_0 \vone \vone^T  \vC^{-1}_0 }{\vone^T  \vC^{-1}_0 \vone}
\right) \vf \\
\nonumber
&\qquad \qquad \times
\left (
\frac{(1 - \Pi[\vc_0(\cdot, X)]\vC_0^{-1} \vone)^2}{\vone^T  \vC^{-1}_0 \vone} \right . \\
\nonumber & \qquad \qquad  +
\Pi\Pi[c_0(\cdot,\cdot)] - \Pi[\vc_0(\cdot, X)]\vC_0^{-1} \Pi[\vc_0(X,\cdot)] 
\biggr ), \\
\label{BCFullStop}
\lefteqn{\Prob\Bigl[\bigl \lvert \hat{\Pi}[f] - \Pi[g] \bigr \rvert \le t_{n-1,0.995} \sqrt{\Var_n[\Pi[g]]}  \Bigr] = 99\%.}
\end{align}
\end{subequations}
Here, $\vone$ is a vector of ones, and $t_{n-1,0.995}$ denotes the $99.5\%$ quantile of the Student-t distribution with $n-1$ degrees of freedom.  For large $n$, $t_{n-1,0.995} \approx 2.58$.  The expressions in \eqref{JHFullBayes} are similar to the conclusion of \BOGOS{}, Proposition 2.  The differences are due to the mean of the Gaussian process being left unspecified in \JH{}, which reduces the degrees of freedom by one, and adds additional terms to the expressions for $\Ex_n[\Pi[g]]$ and $\Var_n[\Pi[g]]$.

Hidden in the definition of $c_0$ is the parameter $\vtheta$. One may place a discrete prior on $\vtheta$, but this strikes us as rather arbitrary.  Thus, in \JH{} we advocate estimating $\vtheta$ by \emph{empirical Bayes}, namely,
\begin{equation} \label{thetaEB}
    \vtheta_{\textup{EB}}
= \argmin_{\vtheta} \biggl \{
\log\left(\vf^T 
\left[ \vC_0^{-1} - 
\frac{ \vC_0^{-1} \vone \vone^T \vC_0^{-1} }{\vone^T \vC_0^{-1} \vone}
\right] \vf 
\right) +  \frac{1}{n} \log(\det(\vC_0))
\biggr \}.
\end{equation}  

\JH{} also presents empirical Bayes as an alternative to assuming the improper prior on $(m,\lambda)$, as discussed in \BOGOS{} Section 4.1.3.  Under empirical Bayes, the posterior marginal for $\Pi[g]$ has the same mean as for the full Bayes approach, but a somewhat smaller variance.  

\JH{} also discusses the alternative of \emph{generalized cross-validation} for estimating the correct covariance function from the integrand data, as alluded to in Section 4.1.2 of \BOGOS{}.  The formulas for the Bayesian cubature and the credible interval width are significantly different than those for full Bayes.

%%%%%%%%%%%%%%%%%%%%%%%%%%%%%%%%%%%%%%
\section{Speeding Up the Computation} \label{sec:Match}
Computing the estimate of $\vtheta$ in \eqref{thetaEB} and then the credible interval according to \eqref{JHFullBayes} involve matrix factorization and computing a matrix determinant, which requires as many as $\Order(n^3)$ operations.  On the other hand, the computational cost of obtaining the integrand data, $\vf$, is $\Order(\$(f)n)$, where $\$(f)$ is the computational cost of  a single integrand value.

If $\$(f)$ is extraordinarily large compared to the expected sample size $n$, then the cost of obtaining integrand data dominates, and the $\Order(n^3)$ cost of matrix operations is unimportant.  However, if $\$(f)$ is close to $\Order(1)$, then the cost of matrix operations may make Bayesian cubature prohibitively costly.

\JH{} presents a scenario where the cost of matrix operations may be reduced to $\Order(n \log n)$ via fast transforms.  The key is choosing covariance functions and cubature nodes that match.  Let the matrix $\vC_0$ be decomposed in terms of its eigenvectors, which comprise the columns of $\vV$, and its eigenvalues, which comprise the diagonal elements of the diagonal matrix $\vLambda$:
\begin{align}
\nonumber
\vC_0 &  = (\vC_1,...,\vC_n) 
= \frac 1n \vV \vLambda \vV^H , 
\quad \vV^H = n \vV^{-1}, \\
\nonumber
\vV &= (\vv_1,...,\vv_n)^T = (\vV_1,...,\vV_n).
\end{align}
Four assumptions are made regarding the covariance function, $c_0$, and the cubature nodes, $\{\vx\}_{i=1}^n$ (\JH{} (25, 27)):
\begin{subequations} \label{fastcompAssump}
	\begin{gather}
	\label{fastcompAssumpA}
	\vV \text{ may be identified analytically}, \\
	\label{fastcompAssumpB}
	\vv_1 = \vV_1 = \vone, \\
	\label{fastcompAssumpC}
	\tilde{\vb}:=\vV^H \vb  \text{ requires only $\Order(n \log(n))$ operations } \forall \vb, \\
	\Pi[c_0(\cdot,\vx)] = 1 \qquad \forall \vx.
	\end{gather}
\end{subequations}
Here, $\vV^H \vb$ is called the \emph{fast transform} of $\vb$ because it takes fewer than the typical $\Order(n^2)$ operations required for matrix-vector multiplication.

\begin{figure}
    \centering
    \includegraphics[height = 4cm]{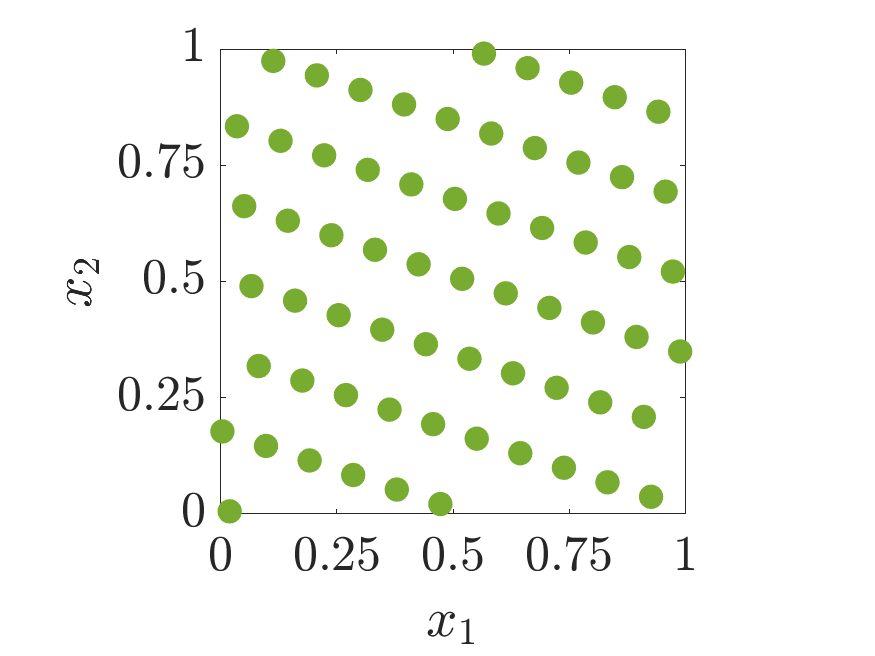} \quad	
    \includegraphics[height = 4.25cm]{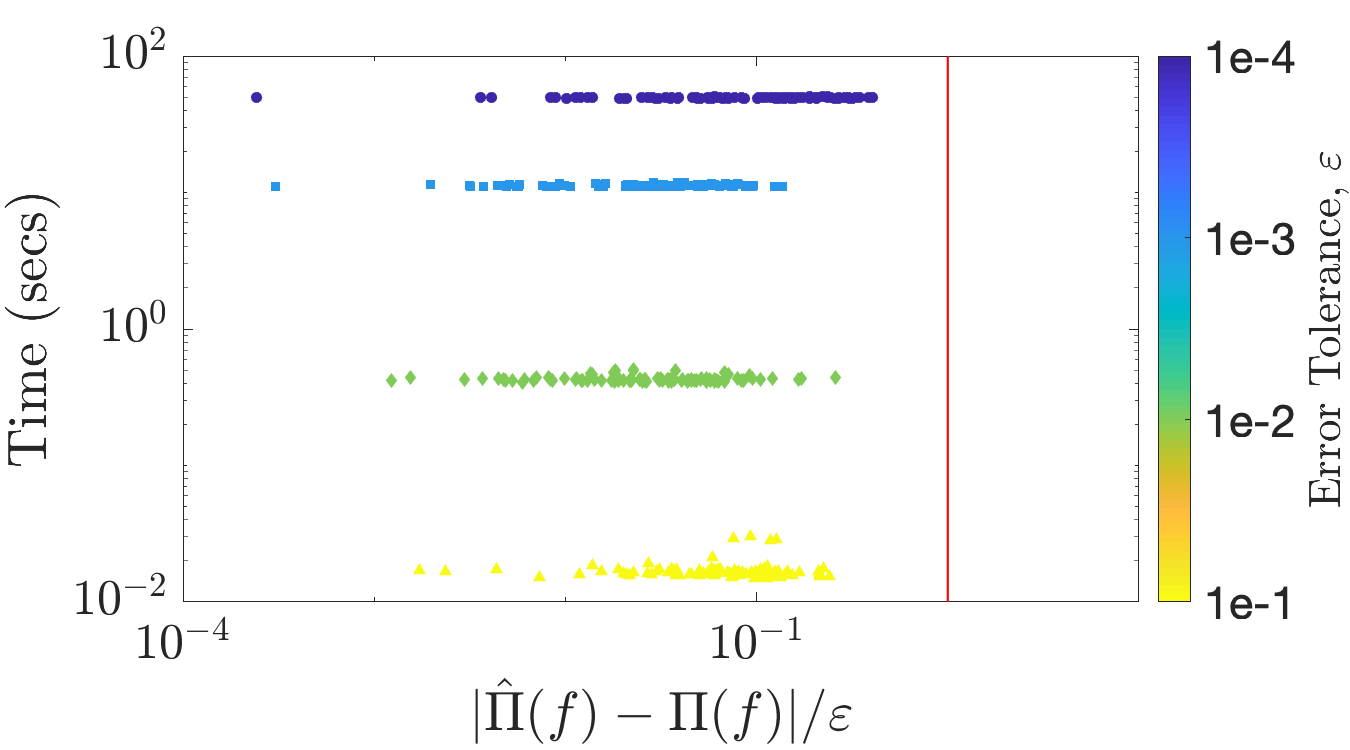}
    \caption{An example of shifted integration lattice nodes in two dimensions (left). The performance of Bayesian cubature for an option pricing example (right).}
    \label{fig:latfig}
\end{figure}

An example of matching covariance functions and cubature nodes is 
\begin{multline*}
    \bullet \text{ Shift-invariant covariance functions, $c_0$, which satisfy}  \\ 
    c_0(\vx, \vx') = \mathring{c}_o(\vx - \vx' \bmod \vone) \qquad \forall \vx, \vx' \in [0,1)^d,
\end{multline*}
for some $\mathring{c}_0$ with period $1$ in each coordinate direction, and
\begin{multline*}
    \bullet \text{ Shifted rank-1 integration lattice node sets, $\{\vx_i\}_{i=1}^n$, which satisfy} \\
    \vx, \vx' , \vx'' \in \{\vx_i\}_{i=1}^n \implies \vx + \vx' - \vx'' \bmod \vone \in \{\vx_i\}_{i=1}^n. 
\end{multline*}
The covariance function in \eqref{the_kernel_eqn_bernoulli} is an example of a shift-invariant covariance function \citep{Hic98b}.  Figure \ref{fig:latfig} (left) depicts a rank-1 integration node set \cite{SloJoe94,DicEtal14a}.  The reason that this family of covariance functions matches this family of cubature nodes and satisfies assumptions \eqref{fastcompAssump} is that the matrix $\vC_0$ is circulant and $\vV$ may be written in terms of complex exponentials.

Under assumptions \eqref{fastcompAssump} one may express \eqref{JHFullBayes} and \eqref{thetaEB} in terms of the fast transforms of the integrand data and the first column of the matrix $\vC_0$ (\JH{} Sections 3.2, 3.3):
\begin{subequations} \label{JHFullBayesFast}
\begin{align}
\label{FastTrans}
\tvf & : = \vV^{H} \vf, \qquad \vell = \diag(\vLambda) = \tilde{\vC}_1 := \vV^{H} \vC_1, \\
    \hat{\Pi}[f] & =  \Ex_n[\Pi[g]] = \frac 1n \sum_{i=1}^n f(\vx_i) = \frac{\tf_1}{n} \qquad \text{the sample average},\\
\Var_n[\Pi[g]] & = \frac{1}{n(n-1)} \left(\frac{\ell_1}{n}  - 1  \right)\sum_{i=2}^n \frac{\bigl \lvert \tf_i\bigr \rvert^2}{\ell_i} , \\
\label{thetaEBfast}
\vtheta_{\textup{EB}} &= 
\argmin_{\vtheta}
\left[
\log\left(
\sum_{i=2}^n \frac{\abs{\widetilde{y}_i}^2}{\ell_i}
\right)  
  + 
\frac{1}{n}\sum_{i=1}^n \log(\ell_i)
\right].
\end{align}
\end{subequations}
Apart from the computations in \eqref{FastTrans}, which require $\Order(n \log n)$ operations, all other calculations in \eqref{JHFullBayesFast} require only $\Order(n)$ operations.  The expression for $\Var_n[\Pi[g]]$ excludes $i=1$ in the sum because we allow $g$ to have an arbitrary constant mean.

Section 5 of \JH{} presents several numerical experiments for Bayesian cubature with shift-invariant covariance functions and lattice nodes sets.  We reproduce one such experiment in Figure \ref{fig:latfig} (right).  The integrand is the payoff of an arithmetic mean Asian option.  The covariance function is the one given by \eqref{the_kernel_eqn_bernoulli} with $r=1$.  A low degree of smoothness is chosen in view of the discontinuities in the partial derivatives of the integrand.  The value of $\gamma$ is determined by empirical Bayes as in \eqref{thetaEBfast}.  The sample size, $n$, is increased in a sequence of powers of $2$ until the stopping criterion implied by \eqref{BCFullStop}, whose terms are computed quickly via \eqref{JHFullBayesFast}, is satisfied.  Four different values of the tolerance were tried, $\varepsilon = 0.1, 0.01, 0.001,$ and $0.0001$. The aim is for the cubature error, $|\Pi[f] - \hat{\Pi}[f]|$, to be no greater than, but not too much less than, the prescribed tolerance nearly all the time.  In this experiment, the error tolerance is always met.  As expected, the computation time increases as the tolerance decreases.  

In this example, and others provided in Section 5 of \JH{}, the cost of evaluating the integrand is modest, and so the cost of obtaining the needed integrand data, $\vf$, is on the same order as the matrix-vector operations required to compute the credible interval.  \JH{} also provides examples of the empirical Bayes and generalized cross-validation approaches to determining the parameters inherent in the covariance function and to using credible intervals as stopping criteria for Bayesian cubature.  All of these approaches are successful, which suggests that they should be explored over a larger range of examples.  

\section{Bayesian Versus Deterministic Analysis} \label{sec:ProbDet}

We return to the situation where the Gaussian process, $g$, has zero mean.  Section 3.2 in \BOGOS{} sets the covariance function, $c$, identical to the reproducing kernel, $k$, of the Hilbert space containing the integrand, $f$, for ``aesthetic'' reasons.  While this makes the application of several results from numerical analysis of deterministic cubature more readily transferable to Bayesian cubature, we think that such a correspondence muddies the waters.  

Suppose that $\{\phi_i\}_{i=1}^\infty$ is an orthonormal basis for the Hilbert space $\ch$ with reproducing kernel $k$.  Then 
\begin{equation}
    \label{kHS}
    k(\vx,\vx') = \sum_{i=1}^\infty \phi_i(\vx) \phi_i(\vx').
\end{equation}
Moreover, the norm of any $f = \sum_{i=1}^n a_i \phi_i \in \ch$ is $\norm[\ch]{f} = \norm[2]{\va}$.  If $g = \sum_{i=1}^\infty A_i \phi_i$, with $A_i \IIDsim \cn(0,1)$, then $g \in \calGP(0,c)$, where $c$ is identical to $k$ as defined in \eqref{kHS}.  While this may seem well and good, note that 
\[
\Ex(\norm[\ch]{g}^2) = \Ex(\norm[2]{\vA}^2) = \norm[2]{(1,1, \ldots) }^2 = \infty.
\]

So setting $c$ identical to $k$ means that we are modeling an integrand with finite norm by a Gaussian process with an infinite expected squared norm.
This seems counter-intuitive.  Nevertheless, there are tantalizing mathematical similarities between the Bayesian and deterministic approaches to cubature.

When the optimal cubature weights are used, the deterministic error bound in \eqref{HJErrBd} may be expressed as
\begin{subequations} \label{HJErrBdTight}
\begin{gather} 
    \label{HJErrbdTightCub}
    \hat{\Pi}[f] = \Pi[\hat{f}] = \Pi[\vk(\cdot, X)] \vK^{-1} \vf, \\
    \label{HJErrBdBdTight}
    \bigl \lvert \hat{\Pi}[f] - \Pi[f] \rvert^2 \le \lVert f - \hat{f}  \rVert_{\calH}^2 \ \bigl \{ \Pi\Pi[k(\cdot,\cdot)] - \Pi[\vk(\cdot, X)]\vK^{-1} \Pi[\vk(X,\cdot)]\bigr \} ,
\end{gather}
\end{subequations}
where $\hat{f}$ is the minimum Hilbert space norm interpolant of the integrand, $f$.  The reason that $\lVert f \rVert_{\calH}$ in \eqref{HJErrBd} can be replaced by $\lVert f - \hat{f} \rVert_{\calH}$ is that the cubature in \eqref{HJErrbdTightCub} integrates $\hat{f}$ exactly. Moreover, $f - \hat{f}$ is orthogonal to $\hat{f}$ under the Hilbert space inner product.  Also note that
\begin{equation} \label{splinenorm}
   \lVert \hat{f}  \rVert^2_{\calH} = \vf^T \vK^{-1} \vf.
\end{equation}

Compare the error bound in \eqref{HJErrBdTight} to Proposition 2 of \BOGOS{}, which implies that 
\begin{subequations} \label{HJBCErrBdTight}
\begin{gather} 
\label{HJBCErrbdTightCub}
    \hat{\Pi}[f] = \Ex[\Pi[g]] = \Pi[\vc_0(\cdot, X)] \vC_0^{-1} \vf , \\
    \Var[\Pi[g]] = \frac{\vf^T \vC_0^{-1} \vf}{n} \ \bigl \{ \Pi\Pi[c_0(\cdot,\cdot)] - \Pi[\vc_0(\cdot, X)]\vC_0^{-1} \Pi[\vc_0(X,\cdot)] \bigr \} , \\
    \label{HJBCErrBdBdTight}
    \Prob\bigl[ \bigl \lvert \hat{\Pi}[f] - \Pi[f] \rvert^2 \le t_{0.995,n}^2 \Var[\Pi[g]] \bigr] = 99\%.
\end{gather}
\end{subequations}
The formulas for the cubature $\hat{\Pi}[f]$ in both the deterministic and Bayesian senses are identical if $k$ is identical to $c_0$.  They are also independent of the vertical scale factor multiplying inherent in the definition of $k$ or $c_0$, i.e., the reproducing kernels $k$ and $47k$ yield the same cubature rule.  Likewise, the bound on the squared error in \eqref{HJErrBdBdTight} and the variance of the integral of the Gaussian process in \eqref{HJBCErrBdBdTight} both contain the common factor $\Pi\Pi[k(\cdot,\cdot)] - \Pi[\vk(\cdot, X)]\vK^{-1} \Pi[\vk(X,\cdot)]$ if $k$ is identical to $c_0$.

Matching the deterministic error bound in \eqref{HJErrBdTight} to the Bayesian credible interval in  \eqref{HJBCErrBdTight} when  $k$ is identical to $c_0$ becomes possible if one applies \eqref{splinenorm} and asserts that 
\begin{equation}
\label{deterasprob}
\lVert f - \hat{f} \rVert_{\calH}^2 \le \frac{ t_{0.995,n}^2 \lVert \hat{f} \rVert^2_{\calH}}{n} = \frac{t_{0.995,n}^2 \vf^T \vK^{-1} \vf}{n}.
\end{equation}
Although this inequality is violated for some $f \in \ch$, it holds for those $f \in \ch$ that are well-modeled by their minimum norm interpolants, $\hat{f}$.  Thus, one can mimic Bayesian cubature via a deterministic cubature which assumes that the integrand satisfies inequality \eqref{deterasprob}.

\section{Further Matters}

A couple of matters deserve further investigation.  How large a family of covariance functions must be considered for effective Bayesian cubature?  A larger family increases the probability that the integrand in question lies in the middle of the space of Gaussian processes used to determine the stopping criterion.  On the other hand, a larger family may require a more tedious choice of the underlying parameters $\vtheta$.

The Bayesian approach to numerical integration assumes a Gaussian process.  Do goodness-of-fit statistics confirm or discredit this assumption?  How does the validity of this assumption affect the reliability of the proposed Bayesian automatic cubature?  The alternative of Student-$t$ processes has been suggested by \cite{ShaWilGha14a}.

Finally, probabilistic numerics---including Bayesian cubature---deserves further participation from statisticians, numerical analysts, and software developers alike.  Statisticians and numerical analysts should become more conversant in each other's language and culture.  Computational problems are better understood when one can look from multiple perspectives.  Moreover, the algorithms that arise from probabilistic numerics should find their way into commonly used software libraries.  Such libraries should be built using  software engineering principles that are familiar to software developers, but perhaps not obvious to statisticians. Our Guaranteed Automatic Integration Library (GAIL) 
\cite{ChoEtal17b} 
is an example of such a library.

\section*{Acknowledgements}
We thank the authors of \BOGOS{} for a provocative article and for urging the computational mathematics community to investigate more carefully probabilistic approaches.  We also thank Sou-Cheng Choi for her comments. This work was supported in part by NSF-DMS-1522687 and NSF-DMS-1638521 (SAMSI).

\bibliographystyle{imsart-nameyear}
\bibliography{FJH23,FJHown23}

\end{document}